\documentclass[11pt]{amsart} 
\usepackage[dvips]{graphics} 
\usepackage{aeguill}
\newtheorem{prop}{Proposition}[section] 
 
\newtheorem{theo}{Theorem}[section] 
\newtheorem{lem}{Lemma}[section] 
\newtheorem{cor}{Corollary}[section] 
\newtheorem{rem}{Remark}[section] 
\newtheorem{defi}{Definition}[section]

\numberwithin{equation}{section}

\newcommand{\R}{\mathbb R}

\begin{document} 
\begin{abstract}
We construct some lift of an almost complex structure to the cotangent bundle, using a connection on the base manifold.  
This unifies the complete lift defined by I.Satô and the horizontal lift introduced  by S.Ishihara and  K.Yano. 
We study some geometric properties of this lift and its compatibility with symplectic forms on the cotangent bundle. 
\end{abstract} 
\title { Almost complex structures on the cotangent bundle } 
\author{Florian Bertrand}
\address{LATP, C.M.I, 39 rue Joliot-Curie 13453 Marseille cedex 13, FRANCE }
\email{bertrand@cmi.univ-mrs.fr}
\maketitle 
\section*{Introduction} 
Analysis on almost complex manifolds recently  became an indispensable tool in 
symplectic geometry with the celebrated work of M.Gromov in \cite{gr}. 
The local existence of pseudoholomorphic discs proved by A.Nijenhuis-W.Woolf 
in their famous paper \cite{ni-wo}, allows to lead some local analysis 
on such manifolds. There is a natural and deep connection beetwen 
local analysis on complex and almost complex manifolds and canonical bundles. 
For instance, the cotangent bundle is tightly related to
extension of biholomorphisms and to the study of stationnary discs. 
Morever, it is well known that the cotangent bundle  plays
a very important  role in symplectic geometry and its applications, since this carries a canonical symplectic structure
induced by the Liouville form. 

Several lifts of an almost complex structure on a base manifold are constructed on the cotangent bundle. These are essentially due to 
I.Satô in \cite{sa} and S.Ishihara-K.Yano in \cite {ya-is}. I.Satô defined a lift 
of the ambient structure as a correction of the
 {\it complete lift};  S.Ishihara-K.Yano introduced the {\it horizontal lift} obtained via a symmetric connection. 
The aim of the present paper is to unifiy and to generalize
these lifts by introducing a more natural almost complex lift called the {\it generalized horizontal lift}. 

It turns out that our construction depends on the introduction of some connection : we study the dependence of the lift on it. 
Our main result states that the structure defined by I.Satô and
 the {\it horizontal lift} are special cases of our general 
construction, obtained by particular choices of  connections (Theorem \ref{propimp}).
We establish some geometric properties of this general lift (Theorems \ref{propprop} and 
\ref{theoholo}). Then we characterize generically the structure constructed by I.Satô  
 by the holomorphicity of the lift of a given diffeomorphism on the bases and by 
the holomorphicity of the complex fiberwise multiplication (Corollary \ref{corocoro} and Corollary \ref{corcor}).



Finally, we study the compatibility between lifted almost complex structures and symplectic forms on the cotangent bundle. 
The conormal bundle of a strictly pseudoconvex hypersurface is a
totally real maximal submanifold in the cotangent bundle endowed with the structure defined by I.Satô. This was
proved by S.Webster (\cite{we}) for the standard complex structure, and by
A.Spiro (\cite{sp}), and  independently by H.Gaussier-A.Sukhov (\cite{ga-su}),   
for the almost complex case. One can search for a symplectic proof of this,  
since every Lagrangian submanifold in a symplectic manifold is totally real for almost complex structures
compatible with the symplectic form.
We prove that for every almost complex manifold 
and every  symplectic form on $T^{*}M$ compatible with the {\it generalized horizontal lift}, 
the conormal bundle of a strictly pseudoconvex hypersurface is not Lagrangian (Proposition \ref{proplag}).  



\section{Preliminaries}
 
Let $M$ be a real smooth manifold of even dimension $n$. We denote by $TM$
and $T^{*}M$ the tangent and cotangent bundles over $M$, by
$\Gamma(TM)$ and $\Gamma(T^{*}M) $ the sets of sections of these
bundles and by $\pi :T^{*}M\longrightarrow M$ the fiberwise
projection. We consider local coordinates systems
$(x_{1},\cdots,x_{n})$ in $M$ and $(x_{1},\cdots,x_{n},p_{1},\cdots,p_{n})$ in $T^{*}M$.
We do not write any sum symbol; we  use Einstein summation convention.
\subsection{Almost complex structures}
\begin{defi} 
An almost complex structure on $M$ is a tensor field $J$ of type
$(1,1)$ which  satisfies $J^{2}=-Id$. The pair $(M,J)$ is called an almost complex manifold.
\end{defi} 
In local coordinates, $J$ is given by $J_{l}^{k}dx^{l}\otimes\partial x_{k}$.

We say that a map $f:(M,J)\longrightarrow (M',J')$ between two
almost complex manifolds is $(J,J')$-holomorphic if :
$$
J'(f(x))\circ d_{x}f=d_{x}f\circ J(x), \mbox{ for every } x \in
M.
$$
If $f:(M,J)\longrightarrow M'$ is a diffeomorphism, we define the
 direct image of $J$ by $f$ by :
$$f_{*}J(y):=d_{f^{-1}(y)}f\circ J(f^{-1}(y))\circ d_{y}f^{-1}, \mbox{ for every } y \in M'.$$
The tensor field $f_{*}J$ is an almost complex structure on $M'$ for which $f$ is $(J,f_{*}J)$-holomorphic.

We recall that the Nijenhuis tensor of the almost complex
structure $J$ is defined by :
$$
N_{J}(X,Y):=[JX,JY]-J[X,JY]-J[JX,Y]-[X,Y] \mbox{ for } X,Y \in
\Gamma (TM).
$$
 

\subsection{Tensors and contractions}

Let $\theta$ be the Liouville form on $T^{*}M$. This one-form is locally
given by $\theta=p_{i}dx^{i}$. The two-form $\omega_{st}:=d\theta$ is the canonical
symplectic form on the cotangent bundle, with local expression
$\omega_{st}=- dx^{k}\wedge dp^{k}.$ We stress out that these forms do not depend on the choice of coordinates on $T^{*}M$.
 
We denote by  $T_{q}^{r}M$ the space of $q$ covariant and $r$ contravariant tensors on $M$.
 For positive $q$, we consider the contraction map
$\gamma : T_{q}^{1}M \rightarrow T_{q-1}^{1}(T^{*}M)$ defined by :
$$
\gamma(R):=
p_{k}R_{i_{1},\cdots,i_{q}}^{k}dx^{i_{1}}\otimes\cdots\otimes
dx^{i_{q-1}}\otimes\partial p_{i_{q}}
$$
for $R=
R_{i_{1},\cdots,i_{q}}^{k}dx^{i_{1}}\otimes\cdots\otimes
dx^{i_{q}}\otimes\partial x_{k}$.
 
 We also define a $q$-form on $T^{*}M$ by  
$\theta (R):=p_{k}R_{i_{1},\cdots,i_{q}}^{k}dx^{i_{1}}\otimes\cdots\otimes dx^{i_{q}}$ for a tensor $R \in  T_{q}^{1}M$ on $M$. 
We notice that $\theta (R)(X_{1},\cdots,X_{q})=\theta(R(d\pi (X_{1}),\cdots,d\pi (X_{q})))$ for
$X_{1},\cdots,X_{q} \in \Gamma(T^{*}M)$.

Since the canonical
symplectic form $\omega_{st}$ establishes a correspondence between
q-forms and  $T_{q-1}^{1}M$, one may define the contraction map $\gamma $ using the Liouville form 
$\theta$ and $\omega_{st}$ by setting, for $X_{1},\cdots,X_{q} \in \Gamma(T^{*}M)$~: 
$$
{}^{t}(\theta (R))(X_{1},\cdots,X_{q})=-\omega_{st}
(X_{1},\gamma (R)(X_{2},\cdots,X_{q})),
$$
where ${}^{t}(\theta (R))(X_{1},\cdots,X_{q})=\theta(R)(X_{2},\cdots,X_{q},X_{1})$.

For a tensor $R \in T_{2}^{1}M$, we have a matricial interpretation
of the contraction $\gamma $; if $R_{i,j}^{k}$ are the coordinates of $R$
then $\gamma(R)$ is given by :
$$\gamma(R)=\left(\begin{array}{ccccc}
 
0 & 0 \\ a^{i}_{j}& 0 \\
  
\end{array}\right) 
\mbox{ } \in \mathcal{M}_{2n}(\R), \mbox{ with
}a^{i}_{j}=p_{k}R_{j,i}^{k}.$$

\subsection{Connections}

Let $\nabla$ be a connection on an almost complex manifold $(M,J)$. We denote by $\Gamma _{i,j}^{k}$ its
Christoffel symbols defined by
$\nabla_{\partial x_{i}}\partial
x_{j}= \Gamma _{i,j}^{k}\partial x_{k}$.
Let also $\Gamma_{i,j}$ defined in local coordinates
$(x_{1}\cdots,x_{n},p_{1},\cdots,p_{n})$ on $T^{*}M$
by the equality $p_{k}\Gamma_{i,j}^{k}=\Gamma_{i,j}$.

The torsion $T$ of $\nabla$ is defined by : 
$$T(X,Y):=\nabla_{X}Y-\nabla_{Y}X-[X,Y],\mbox{ for every } X,Y \in \Gamma
(TM).$$
There are ``natural'' families of connections on an almost complex manifold. 
\begin{defi} A connection $\nabla$ on $M$ is called : 
\begin{enumerate} 
\item almost complex when $\nabla_{X}(JY)=J\nabla_{X}Y$ for every $X,Y$ $\in
\Gamma(TM)$, 
\item minimal when its torsion $T$ is equal to $\frac{1}{4}N_{J}$, 
\item symmetric when its torsion $T$ is identically zero.
\end{enumerate} 

\end{defi}    

A.Lichnerowicz proved, in \cite{li}, that for any almost complex manifold, 
the set of almost complex and minimal connections is nonempty. This fact is crucial in the following.

We introduce a tensor $\nabla J \in T_{2}^{1}M$ which measures the lack of
complexity of the connection $\nabla$ : 
\begin{equation}\label{equation1}
(\nabla J)(X,Y):=\nabla_{X}JY-J\nabla_{X}Y \mbox{ for every } X,Y \in \Gamma(TM).
\end{equation}
Locally we have $(\nabla J)_{i,j}^{k}= \partial x_{i}J^{k}_{j} -
J^{k}_{l}\Gamma_{i,j}^{l} + J^{l}_{j}\Gamma_{i,l}^{k}$.  
 
To the connection $\nabla$ we associate three other connections :
\begin{itemize} 
 
\item $\overline{\nabla}:=\nabla-T$. The Christoffel symbols
$\overline{\Gamma}_{i,j}^{k}$ of
$\overline{\nabla}$ are given by
$\overline{\Gamma}_{i,j}^{k}=\Gamma_{j,i}^{k}.$
 
\item 
$\widetilde{\nabla}:=\nabla-\frac{1}{2}T.$ The connection $\widetilde{\nabla}$
is a symmetric connection and its Christoffel symbols
$\widetilde{\Gamma} _{i,j}^{k}$ are given by : $\widetilde{\Gamma}
_{i,j}^{k}=\frac{1}{2}(\Gamma _{i,j}^{k}+\Gamma _{j,i}^{k}).$
 
\item a connection on $(M,T^{*}M)$, still denoted by $\nabla$, and defined by :
$$(\nabla_{X}s)(Y):=X.s(Y)-s\nabla_{X}Y\mbox{ for every } X,Y \in \Gamma
(TM) \mbox{ and } s \in \Gamma (T^{*}M).$$
 
\end{itemize}


Let $x \in M$ and let $\xi \in T^{*}M$ be such that
$\pi(\xi)=x$. The horizontal distribution $H^{\nabla}$ of $\nabla$ is defined by :
$$H^{\nabla}_{\xi}:=\{d_{x}s(X),\mbox{ }X \in T_{x}M, s \in \Gamma
(T^{*}M), s(x)=\xi, \nabla_{X}s=0\} \subseteq T_{\xi}T^{*}M.$$
We recall that $d_{\xi}\pi$ induces an isomorphism
between $H^{\nabla}_{\xi}$ and $T_{x}M.$ Moreover we have the
following decomposition : $T_{\xi}T^{*}M=H^{\nabla}_{\xi}\oplus
T_{x}^{*}M.$ So an element $Y \in T_{\xi }T^{*}M$ decomposes as
$Y=(X,v^{\nabla}(Y))$, where $v^{\nabla}: T_{\xi }T^{*}M
\longrightarrow T_{x}^{* }M$ is the projection on the vertical space
$T^{*}_{x}M$ parallel to $H^{\nabla}_{\xi }$.

\section{Generalized horizontal lift on the cotangent bundle}

Let $(M,J)$ be an almost complex manifold. We first recall the definitions of the structures constructed by I.Satô 
and S.Ishihara-K.Yano. Then we introduce a new almost complex
lift of $J$ to the cotangent bundle $T^*M$ over $M$ and we prove that this unifies the  complete lift 
 and the horizontal lift.

\subsection{Complete and horizontal lifts}

We consider  the complete lift denoted by $J^{c}$ and defined by 
I.Satô in \cite{sa} as follows :    
let $\theta(J)$ be the one-form on $T^{*}M$ with local expression
$\theta(J)= p_{k}J_{l}^{k}dx^{l}$. We define $J^{c}$ by the identity $d (\theta(J)) =
\omega_{st}(J^{c}.,.).$ Then $J^{c}$ is locally given by :
$$J^{c}=\left(\begin{array}{ccccc}
 
J_{j}^{i} & 0 \\
p_{k}(\partial x_{j}J_{i}^{k}-\partial x_{i}J_{j}^{k})&
J_{i}^{j}\\
\end{array}\right).$$ 
The complete lift $J^{c}$  is an
almost complex structure on $T^*M$ if and only if $J$ is an
integrable structure on $M$, that is if and only if $M$ is a complex
manifold.
Introducing a correction term which involves the non integrability of $J$, I.Satô 
 obtained  an almost complex structure
on the cotangent bundle  (\cite{sa}); this is given by : 
$$
\widetilde{J}:=J^{c}-\frac{1}{2}\gamma(JN_{J}).
$$
For convenience we will also call $\widetilde{J}$  the {\it complete lift} of $J$.  
The coordinates of $JN_{J}$ 
are given by : 
$$
JN_{J}(\partial x_{i},\partial x_{j}) =
[-\partial x_{j}J^{k}_{i}+\partial x_{i}J^{k}_{j} +
J^{k}_{s}J^{q}_{i}\partial x_{q}J^{s}_{j} -
J^{k}_{s}J^{q}_{j}\partial x_{q}J^{s}_{i}]dx^{k}.
$$
Thus we have the following local expression of $\widetilde{J}$ :
$$
\widetilde{J}=\left(\begin{array}{ccccc}
 
J_{j}^{i} & 0 \\ B^{i}_{j}& J_{i}^{j}\\
\end{array}\right),\mbox{ with } B^{i}_{j} =
\frac{p_{k}}{2}[\partial x_{j}J^{k}_{i} -
\partial x_{i}J^{k}_{j} +
J^{k}_{s}J^{q}_{i}\partial x_{q}J^{s}_{j} -
J^{k}_{s}J^{q}_{j}\partial x_{q}J^{s}_{i}].
$$

\vspace{0.5cm}

 We now recall the definition  of the horizontal lift of an almost complex structure. 
Let $\nabla$ be a connection on $M$ and $\widetilde{\nabla}:=\nabla -\frac{1}{2}T$. 
The horizontal lift of $J$ is defined in \cite{ya-is} by :
$$J^{H,\nabla}:=J^{c}+\gamma ([\widetilde{\nabla}J]),$$ 
where the tensor $[\widetilde{\nabla}J] \in T_{2}^{1}M$ is given by :
$$[\widetilde{\nabla}J](X,Y):=-(\widetilde{\nabla}J)(X,Y)+(\widetilde{\nabla}J)(Y,X), 
\mbox{ for every } X,Y \in \Gamma (TM)\mbox{ } (\widetilde{\nabla}J \mbox{ is defined in (\ref{equation1})}).$$ 
 
 S.Ishihara and K.Yano proved that $J^{H,\nabla}$ is an almost
complex structure on $T^{*}M$. It is important to notice that without symmetrizing $\nabla$, the horizontal lift of $J$ is 
not an almost complex structure.
The structure $J^{G,\nabla}$ is locally given by :
$$J^{H,\nabla}=\left(\begin{array}{ccccc}
 
J^{i}_{j}& 0 \\
\widetilde{\Gamma}_{i,l}J^{l}_{j}-
\widetilde{\Gamma}_{j,l}J^{l}_{i}& J^{j}_{i}\\
\end{array}\right).$$

The complete and the horizontal lifts are both a correction of $J^{c}.$ Our aim is 
to unify and to characterize these two almost complex structures.

\subsection{Construction of the generalized horizontal  lift}

Let $x \in M$ and let $\xi \in T^{*}M$ be such that
$\pi(\xi)=x$.
Assume that $H$ is a distribution satisfying the local decomposition $T_{\xi}T^{*}M=H_{\xi}\oplus T_{x}^{*}M$. 
From an algebraic point of view it is natural to lift the almost complex structure $J$ as a product structure, that is  
 $J\oplus {}^{t}J$ with respect to $H_{\xi}\oplus T_{x}^{*}M$. 
Since any such distribution determines and is determined by a unique connection 
one may define a lifted almost complex structure using a connection 
(this point of view is inspired by  P.Gauduchon in ~\cite{gd}).

Let $\nabla$ be a connection on $M$. We consider the connection induced by
$\nabla$ on $(M,T^{*}M)$, defined in subsection~2.3.  For a vector
$Y=(X,v^{\nabla}(Y)) \in T_{\xi }T^{*}M=H^{\nabla}_{\xi}\oplus
T_{x}^{*}M$, we define~:
$$
J^{G,\nabla}(Y)=(JX,{}^{t}J(v^{\nabla}(Y))),
$$
where $JX=(d_{\xi}\pi_{|H^{\nabla}_{\xi }})^{-1}(J(x)d_{\xi }\pi(X))$.
 

\begin{defi} 
The almost complex structure $J^{G,\nabla}$ is called the  generalized horizontal lift 
of $J$ associated to the connection $\nabla$. 
\end{defi}

We first study the dependence of $J^{G,\nabla}$ on the connection $\nabla$.
 
\begin{prop}\label{propgaud} 
Assume that $\nabla$ and $\nabla'$ are two connections on $(M,J)$.
Then $J^{G,\nabla}=J^{G,\nabla'}$ if and only if the tensor
$L:=\nabla'-\nabla$ satisfies $L(J.,.)=L(.,J.)$.
 
\end{prop} 

\begin{proof}.  
Let $\nabla$ and $\nabla'$ be two connections on $(M,J)$ and let $L \in T^{1}_{2}(M)$ be 
the tensor defined by $L:=\nabla'-\nabla$. We notice that, considering the induced
connections on $(M,T^{*}M)$, we have :
$$
\nabla'_{X}s=\nabla_{X}s-s(L(X,.)).
$$
Moreover~:
$$
v^{\nabla'}(Y)=v^{\nabla}(Y)-\xi (L(d_{\xi
}\pi(X),.)),$$
where $Y=(X,v^{\nabla}(Y)) \in T_{\xi }T^{*}M$.

A vector $Y \in T_{\xi}T^{*}M$ can be written $Y=(X,v^{\nabla}(Y))$ in the decomposition
$H^{\nabla}_{\xi } \oplus T_{x}^{*}M$ of $T_{\xi}T^{*}M$ and
$Y=(X',v^{\nabla'}(Y))$ in  $H^{\nabla'}_{\xi}
\oplus T_{x}^{*}M$, with $d_{\xi }\pi(X)=d_{\xi }\pi(X')$. By 
construction we have  $d_{\xi}\pi(JX)=d_{\xi}\pi(JX')$. Thus $J^{G,\nabla'}=J^{G,\nabla}$ if and only
if $v^{\nabla}(J^{G,\nabla'}Y)=v^{\nabla}(J^{G,\nabla}Y)$ for every $\xi \in T^{*}M$ and $Y \in T_{\xi}T^{*}M$. 
Let us compute $v^{\nabla}(J^{G,\nabla'}Y)$ :
$$\begin{array}{lll}
v^{\nabla}(J^{G,\nabla'}Y)&=&v^{\nabla'}(J^{G,\nabla'}Y))+\xi(L(Jd_{\xi
}\pi(X),.))\\
 &=&{}^{t}J(v^{\nabla'}(Y))+\xi(L(Jd_{\xi }\pi(X),.))\\
&=&{}^{t}J(v^{\nabla}(Y))-{}^{t}J\xi
(L(d_{\xi}\pi (X),.))+\xi(L(Jd_{\xi }\pi(X),.))\\
&=&v^{\nabla}(J^{G,\nabla}Y)-\xi
(L(d_{\xi}\pi (X),J.))+\xi(L(Jd_{\xi }\pi(X),.)).
 \end{array}
$$
So $J^{G,\nabla'}=J^{G,\nabla}$ if and only if
$L(d_{\xi}\pi(X),J.)=L(Jd_{\xi }\pi(X),.).$ Since $d_{\xi
}\pi_{|H^{\nabla}_{\xi }}$ is a bijection between $H_{\xi }^{\nabla}$
and $T_{x}M$, we obtain the result.
\end{proof}

\vskip 0,2cm
A consequence of Proposition~\ref{propgaud} is the following Corollary :
\begin{cor}\label{coco1}  
Let $\nabla$ and $\nabla'$ be two minimal almost complex
connections. One has $J^{G,\nabla'}=J^{G,\nabla}$.
\end{cor}

\begin{proof}
Since $\nabla$ and $\nabla'$ have the same torsion, the tensor $L:=\nabla-\nabla'$ is symmetric.  Moreover, since
$\nabla$ and $\nabla'$ are almost complex, we have $L(.,J.)=JL(.,.)$.
Thus $L(J.,.)=JL(.,.)=L(.,J.).$
\end{proof}

We see from  Corollary \ref{coco1} that minimal almost complex connections are ``natural'' connections in almost complex
manifolds, to construct generalized horizontal lifts.

\vspace{0.3cm}

The links between the generalized horizontal lift $J^{G,\nabla}$, the complete lift $\widetilde{J}$, 
and the horizontal lift $J^{H,\nabla}$ are given by the following Theorem : 

\begin{theo}\label{propimp} 
We have :

\begin{enumerate} 
\item $J^{G,\nabla}=\widetilde{J}$ if and only if
$S=-\frac{1}{2}JN_{J}$,  where $S(X,Y)=-(\nabla J)(X,Y)+(\nabla J)(Y,X)+T(JX,Y)-JT(X,Y)$, 

\item $J^{G,\nabla}=J^{H,\nabla}$ if and only if $T(J.,.)=T(.,J.)$ and,

\item For every almost complex and minimal connection, we have $J^{G,\nabla}=\widetilde{J}=J^{H,\nabla}$.

\end{enumerate}
\end{theo}

\subsection{Proof of Theorem \ref{propimp}}
The main idea of the proof is to find a tensorial expression of the generalized horizontal structure  $J^{G,\nabla}$,  
involving $J^{c}$. In that way,  we first  describe locally the horizontal distribution $H^{\nabla}$~:
\begin{lem}\label{lemdist} 
We have $H^{\nabla}_{\xi}=\left\{\left(\begin{array}{ccccc}
 
X \\ \Gamma _{j,k}X^{j}\\
\end{array}\right), X \in  T_{x}M \right\}$ for $\xi \in T^{*}M$  such that
$\pi(\xi)=x$. 
\end{lem}

\begin{proof} 
 Let us prove that $H^{\nabla}_{\xi}\subseteq \left\{\left(\begin{array}{ccccc}
X \\ \Gamma _{j,k}X^{j}\\
\end{array}\right), X \in T_{x}M \right\}.$ 
 Let $Y \in H_{\xi }^{\nabla}$; $Y$ is equal to $d_{x}s(X)$ where $X \in T_{x}M$ and $s$ is a
 section of the cotangent bundle such that $\nabla_{X}s=0$.  
Locally we have  $s=s_{i}dx^{i}$, $X=X^{i}\partial{x_{i}}$ and  so :
 $$Y=\left(\begin{array}{ccccc}
 
X \\ X^{j}\partial{x_{j}}s_{i}\\
\end{array}\right).$$  
Since $\nabla_{X}s=0$ we obtain :
$$0=X^{j}\nabla_{\partial
x_{j}}(s_{i}dx^{i})=X^{j}s_{i}\nabla_{\partial x_{j}}dx^{i}
+X^{j}\partial x_{j}s_{i} dx^{i}=-X^{j}s_{i}\Gamma
_{j,k}^{i}dx^{k}+X^{j}\partial  x_{j}s_{k} dx^{k}.$$ Therefore
$X^{j}\partial x_{j}s_{k}=X^{j}s_{i}\Gamma _{j,k}^{i}=X^{j}\Gamma
_{j,k}$. This proves the inclusion.
 
Moreover the following decomposition insures  the equality :
$$T_{\xi}T^{*}M=\left\{\left(\begin{array}{ccccc}
 
X \\ \Gamma _{j,k}X^{j}\\
\end{array}\right), X \in T_{x}M \right\}\oplus T_{x}^{*}M.$$

\end{proof} 
\skip 0.1cm

\vskip 0.1cm
The following Proposition gives the local expression of the generalized horizontal lift which is necessary
to obtain the desired tensorial expression stated in part (2).

\begin{prop} \label{propgaudu} \mbox{ }
\begin{enumerate} 

\item  With respect to the local coordinates system $(x_{1},\cdots,x_{n},p_{1},\cdots,p_{n})$, $J^{G,\nabla}$ 
is given by :  
$$J^{G,\nabla}=\left(\begin{array}{ccccc} 
J^{i}_{j}& 0 \\
\Gamma_{l,i}J^{l}_{j}-
\Gamma_{j,l}J^{l}_{i}& J^{j}_{i}\\
\end{array}\right).$$

\item We have $J^{G,\nabla}=J^{c}+\gamma(S)$  with
$S(X,Y)=-(\nabla J)(X,Y)+(\nabla J)(Y,X)+T(JX,Y)-JT(X,Y)$.
\end{enumerate}
\end{prop}

\begin{proof} 
We first prove part (1).  We denote by $\delta ^{i}_{j}$ the Kronecker symbol. 
With respect to the local coordinates system $(x_{1},\cdots,x_{n},p_{1},\cdots,p_{n})$, the 
structure $J^{G,\nabla}$ is locally given by  : 
$$J^{G,\nabla}=\left(\begin{array}{ccccc}
J^{i}_{j} & 0 \\ a^{i}_{j} & J^{j}_{i} \\  
\end{array}\right).$$   
Since  $\left(\begin{array}{ccccc} 
\delta^{j}_{i} \\ \Gamma _{i,j}\\
\end{array}\right) \in H^{\nabla}_{\xi}$, 
it follows from Lemma \ref{lemdist}, that for every $i \in
\{1,\cdots,n\}$ : 
$$J^{G,\nabla}\left(\begin{array}{ccccc} 
\delta^{j}_{i} \\ \Gamma _{i,j}\\
\end{array}\right)=\left(\begin{array}{ccccc}   
J^{j}_{i} \\ \Gamma _{k,j}J^{k}_{i}\\
\end{array}\right).$$ 
Hence we have : 
$a^{i}_{j}=\Gamma_{l,i}J^{l}_{j}- \Gamma_{j,l}J^{l}_{i}.$ This concludes the proof of  part (1).

Then we prove part (2). Using  the local expression of $J^{c}$, we have :
$$J^{G,\nabla}=J^{c}+\left(\begin{array}{ccccc} 
0& 0 \\
-p_{k}\partial x_{j}J_{i}^{k}+p_{k}\partial x_{i}J_{j}^{k}+\Gamma_{l,i}J^{l}_{j}-
\Gamma_{j,l}J^{l}_{i}& 0\\
\end{array}\right).$$ 
Since  $\nabla_{\partial x_{i}}(J\partial x_{j}) 
=\partial x_{i}J_{j}^{k}\partial x_{k}+\Gamma_{i,l}^{k}J^{l}_{j}\partial x_{k}$, it follows  that :
 $$-p_{k}\partial x_{j}J_{i}^{k}+p_{k}\partial x_{i}J_{j}^{k}+\Gamma_{l,i}J^{l}_{j}-
\Gamma_{j,l}J^{l}_{i}=
p_{k}dx^{k}[-\nabla_{\partial x_{j}}(J\partial x_{i})+\overline{\nabla}_{\partial x_{i}}(J\partial x_{j})]
.$$
We define $S'(X,Y):=-\nabla_{X}(JY)+\overline{\nabla}_{Y}(JX)=-\nabla_{X}(JY)+\nabla_{Y}{JX}+T(JX,Y)$
 and we notice that
 $S'(\partial x_{i},\partial x_{j})=-\nabla_{\partial x_{i}}(J\partial x_{j})+
\overline{\nabla}_{\partial x_{j}}(J\partial x_{i}).$
We point out that  $S'$ is not a tensor. However  with a correction term, we obtain the tensor $S$ :  
$$\begin{array}{lll}
S(X,Y)&=&S'(X,Y)+J[X,Y]\\
&=& -\nabla_{X}(JY)+\nabla_{Y}(JX)+T(JX,Y)+J\nabla_{X}Y-J\nabla_{Y}X-JT(X,Y)\\ 
&=&-(\nabla J)(X,Y)+(\nabla J)(Y,X)+T(JX,Y)-JT(X,Y).
\end{array}
$$
The components of $S$ are given by $S(\partial x_{i},\partial x_{j})=S'(\partial x_{i},\partial x_{j})$
 and so
$J^{G,\nabla}=J^{c}+\gamma(S)$. 
\end{proof}
 \vspace{0.3cm}

Hence we may compare the three lifted structures via their intrinsic expressions given by :
\begin{itemize} 
\item $J^{G,\nabla}=J^{c}+\gamma (S)$   (Proposition \ref{propgaudu}),
\item $\widetilde{J}=J^{c}-\frac{1}{2}\gamma (JN_{J})$  (see subsection 2.2) and,    
\item $J^{H,\nabla}=J^{c}+\gamma ([\widetilde{\nabla}J])$ with
$[\widetilde{\nabla}J](X,Y)=-(\widetilde{\nabla}J)(X,Y)+(\widetilde{\nabla}J)(Y,X)$ (see subsection 2.2).
\end{itemize} 
The lecture of the two first expressions gives part (1) of Theorem \ref{propimp}. 

To prove  (2), we  notice that : 
$$
\begin{array}{lll}
[\widetilde{\nabla}J](X,Y)&=&-(\widetilde{\nabla}J)(X,Y)+(\widetilde{\nabla}J)(Y,X)\\ 
&=&-(\nabla J)(X,Y)+(\nabla J)(Y,X)+\frac{1}{2}T(X,JY)+\frac{1}{2}T(JX,Y)-JT(X,Y).
\end{array}$$

Let us prove part (3)  of Theorem \ref{propimp}. The equality $J^{G,\nabla}=\widetilde{J}$  
follows from  the fact that $\nabla J=0$ because the connection $\nabla$ is almost complex and from the equality 
$-T(J.,.)+JT(.,)=\frac{1}{4}JN_{J}+\frac{1}{4}JN_{J}=\frac{1}{2}JN_{J}.$
Since $T=\frac{1}{4}N_{J}$ and  $N_{J}(J.,.)=N_{J}(.,J.)$ we have $J^{G,\nabla}=J^{H,\nabla}$.  

The proof of  Theorem \ref{propimp} is now achieved. \qed

\vskip 0.3cm 

We end this section with :  

\begin{cor}\label{lemmeegal} 
We have $J^{H,\nabla}=J^{G,\widetilde{\nabla}}.$
\end{cor} 
  
\begin{proof} 

This is a direct consequence of Theorem \ref{propimp} since $J^{H,\nabla}=J^{H,\widetilde{\nabla}}$ and 
$J^{G,\widetilde{\nabla}}=J^{H,\widetilde{\nabla}}$ by part (2).

\end{proof} 
\skip 0.1cm

We point out that Corollary \ref{lemmeegal}  may also be proved using  Lemma \ref{lemdist} and the distribution $D$ of horizontal lifted vectors defined by 
 S.Ishihara and  K.Yano as follows : let $x \in M$ and $\xi  \in T^{*}M$ such that $\pi(\xi)=x$. Assume  
$X^{H,\nabla}$ is the horizontal lift of $X \in T_{x}M$  on the cotangent bundle given 
in \cite{ya-is} by :
$$ X^{H,\nabla}=\left(\begin{array}{ccccc}
 
X \\ \widetilde{\Gamma} _{j,k}X^{j} \\
\end{array}\right) \in T_{\xi}T^{*}M.$$
Then  the distribution $D$ of horizontal lifted vectors is defined by $D_{\xi}=\{X^{H,\nabla}, X \in T_{x}M\}$.
 S.Ishihara and  K.Yano  proved that   $J^{H,\nabla}=J \oplus {}^{t}J$ in the decomposition 
$T_{\xi}T^{*}M=D_{\xi}\oplus T_{x}^{*}M.$ 
From Lemma \ref{lemdist} we have
$D=H^{\widetilde{\nabla}}$ and finally $J^{H,\nabla}=J \oplus {}^{t}J=J^{G,\widetilde{\nabla}}$ with respect to the decomposition 
$T_{\xi}T^{*}M=D_{\xi}\oplus T^{*}_{x}M=H^{\widetilde{\nabla}}_{\xi}
\oplus T^{*}_{x}M.$

\section{Geometric properties of the generalized horizontal lift}

\subsection{Lift Properties}  

In Theorem \ref{propprop} we state the lift properties of the generalized horizontal lift of an almost complex  
structure. 

\begin{theo} \label{propprop}\mbox{ }
\begin{enumerate}  
\item The projection $\pi : T^{*}M \longrightarrow M$ is
$(J^{G,\nabla},J)$-holomorphic.
\item The zero section $s : M \longrightarrow T^{*}M$ is 
$(J,J^{G,\nabla})$-holomorphic.
\item The lift of a diffeomorphism $f :
(M_{1},J_{1},\nabla_{1}) \longrightarrow (M_{2},J_{2},\nabla_{2})$ 
to the cotangent bundle is $(J_{1}^{G,\nabla_{1}},J_{2}^{G,\nabla_{2}})$-holomorphic  if 
and only if $f$ is  a $(J_{1},J_{2})$-holomorphic map  satisfying $f_{*}S_{1}=S_{2}.$
\end{enumerate} 

\end{theo}

We recall that the lift $\widetilde{f}$ of a diffeomorphism
$f : M_{1} \longrightarrow M_{2}$ to the cotangent bundle is defined by
$\widetilde{f}=(f,{}^{t}(df)^{-1})$ and that the differential
$d\widetilde{f}$ is locally given by : 
$$d\widetilde{f}=\left(\begin{array}{ccccc}
 
df & 0 \\ (*) & {}^{t}(df)^{-1} \\
  
\end{array}\right) 
\mbox{ } \in \mathcal{M}_{2n}(\R),$$
where $(*)$ denotes a $(n \times n)$ block of derivatives of $f$ with respect
to $(x_1,\cdots,x_n)$.
  
\begin{proof}[Proof of Theorem \ref{propprop}]
Parts (1) and (2) are consequences of Proposition \ref{propgaudu} (part (1)).

Let us prove part (3). 
Assume that $f : (M_{1},J_{1},\nabla_{1}) \longrightarrow (M_{2},J_{2},\nabla_{2})$ is a
$(J_{1},J_{2})$-holomorphic diffeomorphism satisfying $\widetilde{f}_{*}S_{1}=S_{2}$ and let  $\widetilde{f}$ be its lift to the
cotangent bundle. According to Proposition \ref{propgaudu}, we have $J^{G,\nabla_{i}}=J^{c}+\gamma (S_{i})$ for $i=1,2$.
 We denote by  $\theta _{i}$ and $\omega _{i,st}$  the
Liouville form and the canonical symplectic form of $T^{*}M_{i}$.
 The invariance by lifted diffeomorphisms of these forms
insure  that $\widetilde{f}_{*}\theta _{1}=\theta _{2}$ and
$\widetilde{f}_{*}\omega _{1,st}=\omega _{2,st}$. We also recall that 
${}^{t}(\theta_{i} (S_{i}))=-\omega_{i,st}(.,\gamma (S_{i}).)$.

 Let us establish the following equality
$\widetilde{f}_{*}(J_{1}^{G,\nabla_{1}})=J_{2}^{G,\nabla_{2}}.$  
The first step consists in proving  that the direct image of
$J^{c}_{1}$ by $\widetilde{f}$ is $J^{c}_{2}$.
 By the nondegeneracy of $\omega _{2,st}$, it is equivalent to
obtain the equality  $ \omega _{2,st}(\widetilde{f}_{*}J^{c}_{1}.,.) =\omega
_{2,st}(J^{c}_{2}.,.)$ :
$$
\begin{array}{lll}
\omega _{2,st}(\widetilde{f}_{*}J^{c}_{1}.,.)&=& \omega
_{2,st}(d\widetilde{f}\circ J^{c}_{1}\circ (d\widetilde{f})^{-1}.,.)\\
&=& \omega _{1,st}(J^{c}_{1}\circ
(d\widetilde{f})^{-1}.,(d\widetilde{f})^{-1})\\
&=& \widetilde{f}_{*}(\omega_{1,st}(J^{c}_{1}.,.))\\
&=&\widetilde{f}_{*}d(\theta_{1} (J_{1})),\\
\mbox{ and, }\omega _{2,st}(J^{c}_{2}.,.)&=&d(\theta_{2}
(J_{2})).\\
\end{array}
$$
So let us prove that the
pull-back of $\theta_{2} (J_{2})$ by $\widetilde{f}$ is $\theta_{1}
(J_{1})$. According  to the local expression of
$d\widetilde{f}$, we have
$\widetilde{f}^{*}(\theta_{2}(J_{2}))=\theta_{2} (J_{2}\circ df)$ and then :
$$\widetilde{f}^{*}(\theta_{2} (J_{2}))
=\theta_{2} (df\circ
J_{1})=(\widetilde{f}^{*}\theta_{2}) (J_{1})=\theta_{1}(J_{1}).$$
Thus we obtain  $\widetilde{f}_{*}d(\theta_{1} (J_{1}))=d(\theta_{2}
(J_{2}))$, that is  $ \widetilde{f}_{*}J^{c}_{1}=J^{c}_{2}$. 

To show the result, we may prove that the direct image of
$\gamma(S_{1})$ by $\widetilde{f}$ is
$\gamma(S_{2})$.  We prove more generally that $f_{*}(S_{1})=S_{2}$ if and only if 
$\widetilde{f}_{*}(\gamma(S_{1}))=\gamma(S_{2})$ which is equivalent to
prove that  $f_{*}(S_{1})=S_{2}$ if and only if  $\omega_{2,st}(.,\widetilde{f}_{*}(\gamma (S_{1})).)
=\omega_{2,st}(.,\gamma (S_{2}).) $. We have :
$$
\begin{array}{lll}
\omega _{2,st}(.,\widetilde{f}_{*}\gamma(S_{1}).)&=& \omega
 _{2,st}(.,d\widetilde{f}\circ \gamma(S_{1})\circ
 (d\widetilde{f})^{-1}.)\\
&=& \omega
_{1,st}((d\widetilde{f})^{-1}.,\gamma(S_{1})\circ
(d\widetilde{f})^{-1}.,)\\
&=&\widetilde{f}_{*}(\omega_{1,st}(.,\gamma(S_{1}).))\\
&=&-\widetilde{f}_{*}({}^{t}\theta_{1} (S_{1})).\end{array}$$
Let us check that $f_{*}(S_{1})=S_{2}$ if and only if $\widetilde{f}_{*}{}^{t}(\theta_{1} (S_{1}))={}^{t}(\theta_{2}
(S_{2}))$. We have : 
 $$\widetilde{f}^{*}(\theta_{2} (S_{2})) = \theta_{2}(S_{2}(df,df)) \mbox{ and }
\theta_{1}(S_{1})= (\widetilde{f}^{*}\theta_{2})(S_{1})=\theta_{2} (df\circ S_{1}).$$
According to this fact and the definition of $\theta(R)$, where $R \in T_{2}^{1}M$ is given in the section $1.2$, 
it follows that  $f_{*}S_{1}=S_{2}$  if and only if $\theta_{2}(S_{2}(df,df))=\theta_{2}(df\circ S_{1})$.
So $f_{*}(S_{1})=S_{2}$ if and only if $\widetilde{f}_{*}(\gamma(S_{1}))=\gamma(S_{2}).$
Finally we have proved that if  $f : (M_{1},J_{1},\nabla_{1}) \longrightarrow (M_{2},J_{2},\nabla_{2})$ is a
$(J_{1},J_{2})$-holomorphic diffeomorphism satisfying $f_{*}S_{1}=S_{2}$  then
$\widetilde{f}$ is $(J_{1}^{G,\nabla_{1}},J_{2}^{G,\nabla_{2}})$-holomorphic.   

Reciprocally if $\widetilde{f}$ is $(J_{1}^{G,\nabla_{1}},J_{2}^{G,\nabla_{2}}$)-holomorphic  then $f$ is 
$(J_{1},J_{2})$-holomorphic. 
Indeed  the zero section $s_{1}:M_{1}\longrightarrow T^{*}M_{1}$ is $(J_{1},J_{1}^{G,\nabla_{1}})$-holomorphic 
by part (2) of Theorem \ref{propprop}, 
the projection  $\pi_{2}:T^{*}M_{2}\longrightarrow M_{2}$ is $(J_{2}^{G,\nabla_{2}},J_{2})$-holomorphic 
by part (1) of Theorem \ref{propprop} and we have the equality $f=\pi_{2}\circ \widetilde{f}\circ s_{1}$. 
Since $f$ is $(J_{1},J_{2})$-holomorphic we have $\widetilde{f}_{*}J_{1}^{c}=J_{2}^{c}$. Then the 
$(J_{1}^{G,\nabla_{1}},J_{2}^{G,\nabla_{2}}$)-holomorphicity of $\widetilde{f}$ implies the equality 
$\widetilde{f}_{*}(\gamma (S_{1}))=\gamma (S_{2})$, that is $f_{*}S_{1}=S_{2}.$    

\end{proof}
\skip 0.1cm

As a Corollary, we obtain the lift properties of the  complete  and the horizontal 
lifts by considering special connections. We point out that Theorem \ref{propprop} 
and Corollary \ref{corocoro} characterize the  complete lift via the lift of diffeomorphisms. 
 
\begin{cor} \label{corocoro}\mbox{ }
\begin{enumerate}  
\item The lift of a diffeomorphism $f :
(M_{1},J_{1}) \longrightarrow (M_{2},J_{2})$ to the cotangent bundle
is $(\widetilde{J_{1}},\widetilde{J_{2}})$-holomorphic if and only if $f$ is  $(J_{1},J_{2})$-holomorphic.
\item  The lift of a diffeomorphism $f :
(M_{1},J_{1},\nabla_{1}) \longrightarrow (M_{2},J_{2},\nabla_{2})$ 
to the cotangent bundle is $(J_{1}^{H,\nabla_{1}},J_{2}^{H,\nabla_{2}})$-holomorphic  if 
and only if $f$ is  a $(J_{1},J_{2})$-holomorphic map  satisfying 
$f_{*}[\widetilde{\nabla_{1}}J_{1}]=[\widetilde{\nabla_{2}}J_{2}].$

\end{enumerate} 

\end{cor}

\begin{proof} 

To prove  part (1), we consider  almost complex and minimal connections $\nabla_{1}$
 and $\nabla_{2}$ on  $M_{1}$ and $M_{2}.$ Hence  $\widetilde{J_{1}}=J^{G,\nabla_{1}}=
J_{1}^{c}+\gamma(S_{1})$ and 
$\widetilde{J_{2}}=J^{G,\nabla_{2}}=J^{c}+\gamma(S_{2}).$ We have $S_{1}=-\frac{1}{2}J_{1}N_{J_{1}}$
 and $S_{2}=-\frac{1}{2}J_{2}N_{J_{2}}.$ 
We notice that if  $f :(M_{1},J_{1}) \longrightarrow (M_{2},J_{2})$  is  a $(J_{1},J_{2})$-holomorphic diffeomorphism then 
$f_{*}N_{J_{1}}=N_{J_{2}}$ and then $f_{*}J_{1}N_{J_{1}}=J_{2}N_{J_{2}}.$ According to Theorem \ref{propprop} 
the lift of a diffeomorphism $f$ to the cotangent bundle is $(\widetilde{J_{1}},\widetilde{J_{2}})$-holomorphic if and only if 
$f$ is $(J_{1},J_{2})$-holomorphic.

Finally, part (2) follows from the equality $J^{G,\widetilde{\nabla}}=J^{H,\nabla}$ obtained in Corollary \ref{lemmeegal} and from
 Theorem \ref{propprop}.  
  
\end{proof}

We point out that the projection (resp. the zero section) is $(J',J)$-holomorphic (resp $(J,J')$-holomorphic)
 for $J'=\widetilde{J},J^{H,\nabla}$ due to local expressions of the 
 complete lift and of the horizontal lift.

\skip 0.1cm

\subsection{Fiberwise multiplication }
We consider the multiplication map $ Z : T^{*}M \longrightarrow T^{*}M $  by
a complex number $a+ib$ with $b\neq0$ on the cotangent bundle. This is locally defined by  $Z(x,p)=(x,(a+b{}^{t}J(x))p)$. 
For $(x,p) \in T^{*}M$ we have
$d_{(x,p)}Z=\left(\begin{array}{ccccc}
 
Id & 0 \\ C & aId+b{}^{t}J\\
\end{array}\right),$ where $C^{i}_{j}=bp_{k}\partial x_{j}J^{k}_{i}.$

\begin{theo}\label{theoholo} 
  The multiplication map $Z$ is $J^{G,\nabla}$-holomorphic if and only if
$(\nabla J)(J.,.)=(\nabla J)(.,J.)$. 
\end{theo} 

\begin{proof} 
Let us evaluate $d_{(x,p)}Z\circ
J^{G,\nabla}(x,p)-J^{G,\nabla}(x,ap+b{}^{t}Jp)\circ d_{(x,p)Z}$. This is equal to : 
$$\left(\begin{array}{ccccc}
 
0 & 0 \\ CJ+(aId+b{}^{t}J)B(x,p)-B(x,ap+{}^{t}Jp)-{}^{t}JC & 0\\
\end{array}\right),$$
$\mbox{where} B^{i}_{j}(x,p)=p_{k}(\Gamma_{l,i}^{k}J^{l}_{j}-\Gamma
_{j,l}^{k}J^{l}_{i}).$

We first notice that
$aB^{i}_{j}(x,p)-B^{i}_{j}(x,ap+b{}^{t}Jp)=
-bp_{k}J^{k}_{s}(\Gamma_{l,i}^{s}J^{l}_{j}-\Gamma
_{j,l}^{s}J^{l}_{i}).$ Let us compute $D=CJ+(aId+b{}^{t}J)B(x,p)-B(x,ap+{}^{t}Jp)-{}^{t}JC$ :  
$$D^{i}_{j}=bp_{k}[\underbrace{J^{l}_{j}\partial x_{l}J^{k}_{i}}_{(1)}+\underbrace{J^{l}_{i}\Gamma
_{s,l}^{k}J^{s}_{j}}_{(2)}-
\underbrace{J^{l}_{i}\Gamma_{j,s}^{k}J^{s}_{l}}_{(2)'}-\underbrace{J^{k}_{s}\Gamma
_{l,i}^{s}J^{l}_{j}}_{(3)}+
\underbrace{J^{k}_{s}\Gamma_{j,l}^{s}J^{l}_{i}}_{(3)'}-\underbrace{J^{l}_{i}\partial x_{j}J^{k}_{l}}_{(1)'}].$$
We obtain $(1)+(2)+(3)=J^{l}_{j}(\partial x_{l}J^{k}_{i}+J^{s}_{i}\Gamma_{l,s}^{k}-J^{k}_{s}\Gamma_{l,i}^{s})$
and
$(1)'+(2)'+(3)'=J^{l}_{i}(\partial x_{j}J^{k}_{l}+J^{s}_{l}\Gamma_{j,s}^{k}-J^{k}_{s}\Gamma_{j,l}^{s}).$
We recognize the coordinates of the tensor $\nabla J$ (section 1.3) :
 $$\partial x_{l}J^{k}_{i}-J^{k}_{s}\Gamma_{l,i}^{s}+J^{s}_{i}\Gamma_{l,s}^{k}=(\nabla J)_{l,i}^{k}
\mbox{ and }
\partial x_{j}J^{k}_{l}-J^{k}_{s}\Gamma_{j,l}^{s}+J^{s}_{l}\Gamma_{j,s}^{k}=(\nabla J)_{j,l}^{k}.$$ 
Finally $D^{i}_{j}=bp_{k}[J^{l}_{j}(\nabla J)_{l,i}^{k}-J^{l}_{i}(\nabla J)_{j,l}^{k}].$
Then $Z$ is $J^{H,\nabla}$-holomorphic if and only if
$J_{j}^{l}(\nabla J)_{l,i}^{k}=(\nabla J)_{j,l}^{k}J^{l}_{i}.$
Since $(\nabla J)_{j,l}^{k}J^{l}_{i}\partial x_{k}=(\nabla J)(\partial x_{j},J\partial x_{i})$
and
$J_{j}^{l}(\nabla J)_{l,i}^{k}\partial x_{k}=(\nabla J)(J\partial x_{j},\partial x_{i})$, 
this concludes the proof of  Theorem \ref{theoholo}.
\end{proof} 
\skip 0.1cm

In particular, the almost complex lift $\widetilde{J}$ may be characterized generically
by the holomorphicity of $Z$; more precisely we have :

\begin{cor}\label{corcor}\mbox{ }   
\begin{enumerate} 
\item The multiplication map $Z$ is $\widetilde{J}$-holomorphic and,
\item $Z$ is $J^{H,\nabla}$-holomorphic if and only if
$(\widetilde{\nabla}J)(J.,.)=(\widetilde{\nabla}J)(.,J.)$.
\end{enumerate} 
\end{cor} 

\begin{proof}
Let us prove part (1). Assume $\nabla$ is an almost complex minimal connection on $M$.
 We have $\widetilde{J}=J^{G,\nabla}$ and by almost complexity of $\nabla$, $\nabla J$ is identically
equal to zero. Theorem \ref{theoholo} implies the $\widetilde{J}$-holomorphicity of $Z$.
 
Part (2) follows from  Theorem \ref{theoholo} and from the equality
$J^{H,\nabla}=J^{G,\widetilde{\nabla}} $ stated in  Corollary 
\ref{lemmeegal}.
 
\end{proof} 

\begin{rem}
In the case of the tangent bundle $TM$, the fiberwise multiplication is holomorphic for the complete lift of $J$ if and only if
  $J$ is integrable. More precisely, ``the lack of holomorphicity'' 
of this map is measured by the Nijenhuis tensor (see \cite{kr}).    
\end{rem} 
\skip 0.1cm

\section{Compatible lifted structures and symplectic forms}

Assume $(M,J)$ is an almost complex manifold. Let $\Gamma=\{\rho=0\}$ be a real smooth hypersurface of $M$, where 
$\rho : M\rightarrow \R$ is a defining function of $\Gamma$.

\begin{defi} \mbox{ }

\begin{enumerate}
\item Let $x \in \Gamma$. The Levi form of $\Gamma $ at $x$ is defined by 
$\mathcal{L}^{J}_{x}(\Gamma)(X)=-d(J^{*}d\rho)(X,JX)$ for any $X \in T_{x}\Gamma$.
\item The hypersurface $\Gamma=\{\rho=0\}$ is strictly $J$-pseudoconvex if its Levi 
form is positive definite at any point $x \in \Gamma$.
\end{enumerate}
\end{defi} 
Let $x \in \Gamma$, we define   
$N^{*}_{x}(\Gamma):=\{p_{x} \in T^{*}_{x}M, (p_{x})_{|T_{x}\Gamma}=0\}$. The {\it conormal bundle} over $\Gamma$,  
defined by the disjoint union $N^{*}(\Gamma):=\bigcup_{x \in \Gamma} N^{*}_{x}(\Gamma)$, is a totally real
 submanifold of $T^{*}M$ endowed with the complete lift (see \cite{ga-su} and  \cite{sp}), that is 
$TN^{*}(\Gamma) \cap \widetilde{J}(TN^{*}(\Gamma))=\{0\}.$  
 To look for a symplectic proof of that fact, 
we search for a symplectic form, $\omega'$, compatible with the complete lift for which $N^{*}(\Gamma )$ 
is Lagrangian, that is $\omega '(X,Y)=0$ for every sections $X,Y$ of $TN^{*}(\Gamma)$.
  More generally we are interested in the compatibility with the 
generalized horizontal  lift. Proposition \ref{proplag} states  that one cannot find such a form.

\begin{prop}\label{proplag}
  Assume  $(M,J,\nabla)$ is an almost complex manifold equipped with a connection. Let  $\omega$ be a
 symplectic form on $T^{*}M$ compatible with the generalized horizontal lift $J^{G,\nabla}$.
There is no  strictly pseudoconvex hypersurface in $M$ whose conormal bundle is Lagrangian with respect to $\omega$.

\end{prop}

\begin{proof}

Let  $\Gamma$  be  a  strictly pseudoconvex  hypersurface in $M$ and let $x \in \Gamma.$ 
Since the problem is purely local we can suppose that  $M=\R^{2m}$, $J=J_{st}+O|(x_{1},\cdots,x_{2m})|$ and $x=0$. 
Since $\Gamma$ is strictly pseudoconvex we can also suppose that $T_{0}\Gamma=\{X \in \mathbb{R}^{2m}, X_{1}=0\}$.
The two-form $\omega$ is  given by $\omega=\alpha_{i,j}dx^{i}\wedge 
dx^{j}+\beta_{i,j}dp^{i}\wedge dp^{j}+\gamma_{i,j}dx^{i}\wedge dp^{j}$. 

Assume that $\omega(X,Y)=0$ for every $X,Y \in TN^{*}(\Gamma )$.  
We have $N^{*}_{0}(\Gamma)=\{p_{0} 
\in T^{*}_{0}\mathbb{R}^{2m}, (p_{0})_{|T_{0}\Gamma}=0\}
=\{(P_{1},0,\cdots,0), P_{1} \in \mathbb{R}\}$. 
Then a  vector $Y \in T_{0}N^{*}(\Gamma)$ can be written
$Y=X_{2}\partial x_{2}+\cdots+X_{2m}\partial x_{2m}+P_{1}\partial p_{1}.$
So we have for $2\leq i<j\leq 2m$ : 
$$
\omega_{(0)}(\partial x_{i},\partial x_{j})=\alpha_{i,j}=0.
$$
Then $w'_{(0)}$ is  given by  
$\omega_{(0)}=\alpha_{1,j}dx^{1}\wedge dx^{j}+\beta_{i,j}dp^{i}\wedge 
dp^{j}+\gamma_{i,j}dx^{i}\wedge dp^{j}.$

Since $J^{G,\nabla}_{(0)}=\left(\begin{array}{ccccc}
J_{st}& 0 \\
 0 & J_{st}\\
\end{array}\right)$ we have $J^{G,\nabla}_{(0)}Y'=\partial x_{2m}$ for 
 $Y'=\partial x_{2m-1}\neq 0 \in T_{0}(T^{*}\Gamma)$. 
Thus $\omega_{(0)}(Y',J^{G,\nabla}_{(0)}Y')=0$ and so 
 $\omega$ is not compatible with $J^{G,\nabla}.$ 

\end{proof}
 Proposition \ref{proplag} is also established for complete and horizontal lifts because
 $J^{G,\nabla}_{(0)}=\widetilde{J}_{(0)}=J^{H,\nabla}_{(0)}$. 
\begin{rem}
 Since the conormal bundle of a (strictly pseudoconvex) hypersurface is Lagrangian for the symplectic form $\omega_{st}$ 
on $T^{*}M$, Proposition \ref{proplag} shows  that $\omega_{st}$ and $J^{G,\nabla}$ are not compatible. 
\end{rem}

\skip 0.1cm

\end{document}